\documentclass[twoside,12pt]{article}

\usepackage{amssymb,amsmath}

\newenvironment{proof}{\begin{trivlist}\item[]{\it
Proof.}}{\hfill$\square$\end{trivlist}}

\newenvironment{proofof}[1]{\noindent{\it Proof of
#1.}}{\hfill$\square$\\\mbox{}}

\newtheorem{theorem}{Theorem}[section]
\newtheorem{corollary}[theorem]{Corollary}
\newtheorem{definition}[theorem]{Definition}
\newtheorem{lemma}[theorem]{Lemma}
\newtheorem{proposition}[theorem]{Proposition}
\newtheorem{remark}[theorem]{Remark}

\def\ca{{\cal A}}

\def\mc{{\mathbb C}}
\def\mz{{\mathbb Z}}
\def\mn{{\mathbb N}_0}
\def\si{\sigma_i}
\def\ti{\tau_i}

\def\detq{\mbox{det}_q}
\def\glnc{GL(N,{\mathbb C})}
\def\mnc{M(N, {\mathbb C})}
\def\omnc{\mathcal{O}(M(N,\mathbb{C}))}
\def\od{{\cal O}(D)}
\def\ok{{\cal O}(K)}
\def\oc{{\cal O}(C)}
\def\ogq{{\cal O}(GL_q)}
\def\ohq{{\cal O}(SL_q)}
\def\omq{{\cal O}(M_q)}
\def\odpergq{{\cal O}(D\setminus GL_q)}
\def\okperhq{{\cal O}(K\setminus SL_q)}
\def\os{{\cal O}(S^2_{q,\infty})}
\def\bx{\beta^{\xi}}
\def\ax{\alpha^{\xi}}
\def\nx{\nu^{\xi}}
\def\psx{\psi^{\xi}}
\def\phx{\varphi^{\xi}}
\def\og{{\cal O}(G)}

\def\om{{\cal O}(M)}
\def\ox{{\cal O}(X)}
\def\oh{{\cal O}(H)}
\def\evx{{\mathrm{ev}}_{\xi}}
\def\pd{\pi_D} 
\def\pc{\pi_C}
\def\pk{\pi_K}
\def\ld{\lambda_D}
\def\lk{\lambda_K}
\def\rk{\rho_K}
\def\im{{\rm im}}
\def\id{\mathrm{id}}

\def\lem{^{\leq m}}

\def\xmin{x_{-1}}
\def\xnul{x_0}
\def\xplus{x_{+1}}

\def\coa{\alpha-GL_q}
\def\cob{\beta-GL_q}

\def\oamq{{\cal O}_{A}(M_q)}
\def\oagq{{\cal O}_A(GL_q)}
\def\oad{{\cal O}_A(D)}
\def\ba{\beta_A}
\def\bax{\beta^{\xi}_A}
\def\krq{K(q)}
\def\klq{K[q,q^{-1}]}
\def\okm{{\cal O}_K(M)}
\def\okg{{\cal O}_K(GL)}
\def\oklqg{{\cal O}_{\klq}(GL_q)} 
\def\okrqg{{\cal O}_{\krq}(GL_q)}
\def\okrqm{{\cal O}_{\krq}(M_q)}
\def\oklqm{{\cal O}_{\klq}(M_q)}
\def\okdperg{{\cal O}_K(D\setminus GL)}
\def\okrqdperg{{\cal O}_{\krq}(D\setminus GL_q)} 
\def\oadperg{{\cal O}_A(D\setminus GL_q)} 
\def\ogl{{\cal O}(GL)}

\date{}
\begin{document}
\title{Orbits for the adjoint coaction on quantum matrices} 
\author{
M. Domokos\!\!
\thanks{This research was supported through a European Community 
Marie Curie Fellowship. Partially supported by OTKA No. F 32325 and 
T 34530.} , \ 
R. Fioresi\!\!
\thanks{Investigation supported by the University of Bologna, funds for 
selected research topics.} , \ and 
T. H. Lenagan}
\maketitle
\begin{abstract} 
Conjugation coactions of the quantum general linear group on the algebra 
of quantum matrices have been introduced in an earlier paper and the 
coinvariants have been determined. In this paper the notion of orbit is 
considered via co-orbit maps associated with $\mc$-points of the space of 
quantum matrices, mapping the coordinate ring of quantum matrices into 
the coordinate ring of the quantum general linear group. 
The co-orbit maps are 
calculated explicitly for $2\times 2$ quantum matrices. 
For quantum matrices of arbitrary size, 
it is shown that when the deformation parameter is transcendental  
over the base field, then the kernel of the co-orbit map associated with 
a $\mc$-point $\xi$ is a right ideal generated by coinvariants, provided that 
the classical adjoint orbit of $\xi$ is maximal. 
If $\xi$ is diagonal with pairwise different eigenvalues, then 
the image of the co-orbit map coincides with the subalgebra of coinvariants 
with respect to the left coaction of the diagonal quantum subgroup of the 
quantum general linear group. 

\end{abstract}




\maketitle


\section{Introduction}\label{intro}

Consider $\omq$, the coordinate ring of $N\times N$ quantum matrices 
over $\mc$, where $q$ is a non-zero element of $\mc$.  
Denote by 
$\ogq$ the coordinate ring of the quantum $\glnc$ 
(see \cite{rtf}). 
Define the map 
\begin{equation}\label{eq:beta}
\beta: \ogq\to\ogq\otimes\ogq,\ \ \ 
\beta(h)=\sum h_2\otimes S(h_1)h_3, 
\end{equation}
where we use Sweedler's convention for the Hopf algebra $\ogq$ (and $S$ 
denotes the antipode). 

Similarly, set 
\[\alpha:\ogq \to\ogq\otimes\ogq,\ \ \ \alpha(h)=\sum h_2\otimes h_3S(h_1).\]  

Both $\alpha$ and $\beta$ are right coactions of the Hopf algebra 
$\ogq$ on itself. 
Being the formal dual of the right adjoint action, $\beta$ is 
usually called the right {\it adjoint coaction} of $\ogq$. 
Obviously $\omq$ is a subcomodule of $\ogq$ with respect to 
$\alpha$ and $\beta$, and we use the same symbols to denote the 
restrictions 
$\alpha: \omq\to\omq\otimes\ogq$, 
$\beta:\omq\to\omq\otimes\ogq$.  
These coactions can be viewed as quantum analogues 
of the adjoint action of $\glnc$ on its Lie algebra $\mnc$. 
The coinvariants with respect to these coactions were described in \cite{dl}. 

Our aim here is to find a counterpart in the quantum setting of the notion of 
orbits of the classical adjoint action. We summarize the results 
concerning $\beta$, the case of $\alpha$ being similar. 
There is a natural way to associate with a 
``$\mc$-point $\xi$ of the space of $N\times N$ quantum matrices'' 
a morphism 
$\bx:\omq\to\ogq$ of right $\ogq$-comodules (the right coaction on $\omq$ is 
given by $\beta$, whereas the right coaction on $\ogq$ comes from the 
comultiplication of $\ogq$). 
When $q$ is specialized to $1$, 
the kernel of $\bx$ is the vanishing ideal of the closure of the 
orbit of $\xi$, so the image of $\bx$ can be identified with the 
coordinate ring of the closure of the orbit of $\xi$. 
The $q$-deformed map $\bx$ is not 
an algebra homomorphism
since $\beta$ is not an algebra homomorphism (see \cite{dl}). 
However, it turns out that 
the kernel of $\bx$ contains the right ideal of $\omq$ generated by the 
elements $\tau_i-\tau_i(\xi)$, $i=1,\ldots,N$, 
where $\tau_1$,\ldots,$\tau_N$ are the basic 
coinvariants for $\beta$ introduced in \cite{dl}. 
When $\xi$ is diagonal, the image of $\bx$ is 
contained in the quantum quotient space $\odpergq$, 
the subalgebra of coinvariants of the left coaction of the diagonal quantum 
subgroup $\od$ on $\ogq$ (note that in the classical case $q=1$, $D$ is the 
stabilizer of $\xi$, provided that $\xi$ has pairwise different diagonal 
entries). 
In the special case when $N=2$ and $q$ is not a root of unity, 
we are able to show that although $\bx$ is not an algebra homomorphism, 
it has further nice properties. 
Namely, if $\xi$ is diagonal, and $\xi$ is not a scalar multiple of any of the 
non-negative powers of $\mathrm{diag}(q^2,1)$, 
then $\bx$ is surjective onto  
$\odpergq$.  
Moreover, for such $\xi$, 
the kernel of $\bx$ is the right ideal of $\omq$ generated 
by the elements $\tau_1-\tau_1(\xi)$, $\tau_2-\tau_2(\xi)$, 
where $\tau_1,\tau_2$ are the basic $\beta$-coinvariants from \cite{dl};
that is, $\tau_1=q^{-2}x_{11}+q^{-4}x_{22}$, a weighted trace, 
and $\tau_2=x_{11}x_{22}-qx_{12}x_{21}$, 
the quantum determinant. When $\xi$ is a scalar multiple of some non-negative 
power of $\mathrm{diag}(q^2,1)$, then the kernel of $\bx$ is larger, and 
the image is finite dimensional. 
In the final section we change the setup, and treat $q$ as an indeterminate 
over a subfield of $\mc$. In this generic situation, both the kernel 
and the image of the co-orbit map coincide with the subset predicted by the 
classical theory, for a $\xi\in M(N,\mc)$ whose classical adjoint orbit 
is maximal. To be more precise, the kernel of the co-orbit map 
is a right ideal generated by coinvariants, and if $\xi$ is diagonal, then the 
image of the co-orbit map is the subalgebra of coinvariants with respect to 
the left coaction of the diagonal quantum subgroup.  
The present work can be viewed as a quantum version of some of the basic 
results of Kostant \cite{k} on the classical adjoint orbits 
of the general linear group acting on its Lie algebra 
(which we use in our proof). 

Recall that $\omq$ is the 
$\mc$-algebra generated by $N^2$ indeterminates $x_{ij}$, for 
$i=1,\dots,N$, subject to the following relations. 
\begin{eqnarray} \label{eq:defining-relations}
x_{ij}x_{il} &=& qx_{il}x_{ij}, \\
x_{ij}x_{kj} &=& qx_{kj}x_{ij}, \notag\\
x_{il}x_{kj} &=& x_{kj}x_{il},  \notag\\
x_{ij}x_{kl} - x_{kl}x_{ij} &=& (q-q^{-1})x_{il}x_{kj},\notag
\end{eqnarray}
for $1 \le i < k \le N$ and $1 \le j < l \le N$, where $q\in\mc^*$. 
The algebra $\omq$ is an iterated Ore extension, and so a noetherian domain. 
The {\em quantum determinant}, $\detq$, is the element 
\[
\detq := \sum_{\sigma \in S_{N}}(-q)^{l(\sigma)}
x_{1,\sigma(1)} \dots x_{N,\sigma(N)}.
\]
It is known that $\detq$ is a central element in $\omq$ 
(see \cite[Theorem 4.6.1]{pw}), and by adjoining its 
inverse we get the {\em quantum general linear group} 
\[
\ogq := \omq[\detq^{-1}].
\]
The algebra $\ogq$ is a Hopf algebra with {\em comultiplication} $\Delta$ 
and  {\em antipode} $S$. Recall that 
$\Delta(x_{ij})=\sum_{k=1}^Nx_{ik}\otimes x_{kj}$. 
We shall need the explicit form of $S$ only in the special case $N=2$, 
when 
$S(x_{11})=x_{22}/\detq$, 
$S(x_{22})=x_{11}/\detq$, 
$S(x_{12})=-q^{-1}x_{12}/\detq$, 
$S(x_{21})=-qx_{21}/\detq$. 
The reader should be aware that in many papers the r\^oles of $q$ and 
$q^{-1}$ are interchanged, and so one has to be careful in translating
results from one paper to another. 
In the special case $q=1$, the algebra $\omq$ becomes the coordinate ring 
$\om$ of $\mnc$, and $\ogq$ becomes the coordinate ring $\ogl$ of $\glnc$. 


\section{Orbits of classical points}\label{orbits} 

Start from the classical situation when we are given a (right) action 
\[X\times G\to X,\ \ \ (x,g)\mapsto xg\] 
of an affine algebraic group $G$ on an affine algebraic variety $X$. 
The orbit of some $x\in X$ is the image of the composition of 
the morphisms 
\[
G\to X\times G\to X,\ \ \ 
g\mapsto(x,g)\mapsto xg. 
\]
Passing to coordinate rings, the action is encoded in the comorphism 
$\mu:\ox\to\ox\otimes\og$, which makes $\ox$ a right $\og$-comodule algebra. 
The point $x$ corresponds to the  $\mc$-algebra homomorphism 
$\mathrm{ev}_x:\ox\to\mc$,  $f\mapsto f(x)$. 
The comorphism of the orbit map 
$G\to xG\subset X$, $g\mapsto xg$ is the composition 
\[\mu_x:
\ox\stackrel{\mu}\longrightarrow\ox\otimes\og\stackrel
{\mathrm{ev}_x\otimes\id}
\longrightarrow\og.
\]
The kernel of $\mu_x$ is the ideal of the (Zariski) closure of the orbit 
$xG$, and so the image of $\mu_x$ can be identified with 
the coordinate ring of the closure of $xG$. 
Note that $\mu_x:\ox\to\og$ is an algebra homomorphism, as well as 
a morphism of right $\og$-comodules (where the right coaction of $\og$ on 
itself comes from the action of $G$ on itself by right translations). 

Now assume in addition that the orbit of $x$ is closed in the Zariski 
topology of $X$, and denote by $H$ the stabilizer of $x$. 
Then $H$ acts by left translation on $G$, and the orbit map 
$G\to xG$ factors through the quotient variety $H\setminus G$, 
and induces an isomorphism of the $G$-varieties 
$H\setminus G\cong xG$. In terms of coordinate rings this means that 
$\im(\mu_x)$ is the subalgebra of $\oh$-coinvariants in $\og$ 
(the action of $H$ on $G$ by left translations induces a left coaction of 
$\oh$ on $\og$). 

Furthermore, if $G$ is reductive, then $\ox^G$, the algebra of 
polynomial invariants (which coincides with the algebra of 
$\mu$-coinvariants) is finitely generated as an algebra, say 
$f_1,\ldots,f_s$ is a generating system. Assume in addition that 
the orbit $xG$ is closed and maximal; that is, $x$ is not contained in 
the closure of another orbit (note that for the adjoint action of 
$\glnc$ on its Lie algebra $\mnc$ this holds for the general $x$).  
Then $\ker(\mu_x)$ can be described in 
terms of polynomial invariants. In this case 
the common zero locus of $f_1-f_1(x),\ldots,f_s-f_s(x)$ is the orbit 
$xG$, so $\ker(\mu_x)$ is the radical 
of the ideal generated by $f_1-f_1(x),\ldots,f_s-f_s(x)$. 

Recall that 
\[\beta=(\id\otimes p)\circ(\tau\otimes \id)\circ
(S\otimes\id\otimes\id)\circ\Delta^{(2)}:
\ogq\to\ogq\otimes\ogq,\]
where 
$\Delta^{(2)}=(\Delta\otimes\id)\circ\Delta=
(\id\otimes\Delta)\circ\Delta$, 
$\tau(h\otimes k)=k\otimes h$, and 
$p:\ogq\otimes\ogq\to\ogq$ is the multiplication map. 
The corresponding formula for $\alpha$ is 
\[\alpha=(\id\otimes p)\circ\tau_{(132)}\circ
(S\otimes\id\otimes\id)\circ\Delta^{(2)},\]
where $\tau_{(132)}(h\otimes k\otimes l)=k\otimes l\otimes h$. 
The above discussion on algebraic group actions 
motivates the following definition. 
Take a surjective $\mc$-algebra homomorphism $\evx:\omq\to\mc$. 
The $\xi$ in the notation refers to the $N\times N$ matrix with complex 
entries obtained by evaluating the given homomorphism on the generators 
$x_{ij}$ of $\omq$. Note that for $\xi\in\mnc$ there is a corresponding 
homomorphism $\evx$ if and only if $\xi$ is a ``quantum matrix'' 
in the sense of \cite{m}; that is,  
if the entries of $\xi$ satisfy the relations defining $\omq$.   
We say that $\xi$ is a {\it $\mc$-point of $M_q$} in this case. 
If $q\neq 1$, then $\xi$ is a $\mc$-point of $M_q$ if and only if there is at 
most one non-zero 
entry in each column and each row of $\xi$, and 
$\xi_{ij}\xi_{kl}\neq 0$ with $i<k$ implies that $j<l$. 
For example, a diagonal $\xi$ satisfies this condition. 

\begin{definition}\label{co-orbit-map} 
Let $\xi$ be a $\mc$-point of $M_q$. 
The {\it co-orbit map} $\bx$ of $\xi$ 
with respect to the right 
coaction $\beta$ is the composition $(\evx\otimes\id)\circ\beta$, that is, 
\[\bx:
\omq\stackrel{\beta}\longrightarrow\omq\otimes\ogq\stackrel{\evx\otimes\id}
\longrightarrow\ogq.\]
In the same manner, $\ax$ is defined to be $(\evx\otimes\id)\circ\alpha$. 
\end{definition}

The above discussion suggests that the kernel of the co-orbit map 
should contain information about the ``embedding of the orbit into $M_q$'', 
and the image of the co-orbit map should determine the isomorphism type 
of (the closure of) the ``orbit'' as a quantum space. 
The main point of this paper is that we treat the co-orbit map 
(and not only its image) as our central object, and demonstrate that 
it has certain nice properties, even though we started 
with a coaction which was not an algebra homomorphism. 
In particular, the coinvariants are used to study the kernel of the 
co-orbit map. 

\begin{remark}\label{co-orbit-in-general} 
{\rm Note that Definition~\ref{co-orbit-map} makes sense if $\beta$ is 
replaced by any right coaction $\nu$ of a Hopf algebra $\ca$ on an 
associative $\mc$-algebra $R$, and $\evx$ is an algebra 
homomorphism of $R$ to $\mc$. 
Proposition~\ref{morphism-of-comodules} below clearly holds in this general 
setting. 
The map $\nu^{\xi}$ (and Proposition~\ref{morphism-of-comodules}) 
appears in \cite{b}, \cite{dk} in the situation when 
$R$ is an $\ca$- comodule algebra. In this case 
$\nx$ is an algebra 
homomorphism as well, and 
$\im(\nx)$ is a right coideal subalgebra of $\ca$. 
A right coideal subalgebra of a quantum group is called a 
{\it quantum homogeneous space}, see \cite[11.6.1]{ks}. 
So if $\ca$ is a quantum group and $R$ is a quantum $\ca$-space 
(that is, $R$ is an $\ca$-comodule algebra), then 
$\im(\nx)$ is a quantum homogeneous space. Moreover, any 
quantum homogeneous $\ca$-space can be obtained as the image of 
$\nx$ for some classical point $\xi$ in some quantum $\ca$-space, 
see \cite[Proposition 1.1]{dk} or \cite[Proposition 3.2]{b}. 
}
\end{remark} 

\begin{proposition}\label{morphism-of-comodules} 
The co-orbit map $\bx:\omq\to\ogq$ is a morphism of right $\ogq$-comodules  
(given by $\beta$ and $\Delta$, respectively). 
Similarly, $\ax$ is a morphism of $\ogq$-comodules (given by $\alpha$ and 
$\Delta$, respectively). 
\end{proposition} 

\begin{proof} The claim is the equality  
$\Delta\circ\bx=(\bx\otimes\id)\circ\beta$. 
Observe that 
$\Delta\circ(\evx\otimes\id)
=(\evx\otimes\id\otimes\id)\circ(\id\otimes\Delta)$, 
because both are equal to 
$\evx\otimes\Delta:\omq\otimes\ogq\to\ogq\otimes\ogq$. 
Since $\beta$ is a right coaction we have 
$(\id\otimes\Delta)\circ\beta=(\beta\otimes\id)\circ\beta$. 
Using these two equalities one gets 
\begin{align*}
\Delta\circ\bx
&=\Delta\circ(\evx\otimes\id)\circ\beta
\\&=(\evx\otimes\id\otimes\id)\circ(\id\otimes\Delta)\circ\beta  
\\&=(\evx\otimes\id\otimes\id)\circ(\beta\otimes\id)\circ\beta
\\&=(\bx\otimes\id)\circ\beta.
\end{align*}  
\end{proof} 
\goodbreak

The spaces of coinvariants 
\begin{align*}
\omq^{\coa}:=\{f\in\omq\mid \alpha(f)=f\otimes 1\} 
\\ \omq^{\cob}:=\{f\in\omq\mid \beta(f)=f\otimes 1\}
\end{align*} 
are studied in \cite{dl}. 
It is shown there that although $\alpha$ and $\beta$ are not algebra 
homomorphisms, nevertheless 
\begin{equation}\label{alpha-of-product} 
\alpha(fh)=\alpha(f)\alpha(h)\mbox{ \ if \ }h\in\omq^{\coa}
\end{equation} 
and similarly 
\begin{equation}\label{beta-of-product} 
\beta(fh)=\beta(f)\beta(h)\mbox{ \ if \ }f\in\omq^{\cob}. 
\end{equation} 
It follows that $\omq^{\coa}$ and $\omq^{\cob}$ are subalgebras. 
Let us recall their generators. 
Fix an integer $t$ with $1\le t \le N$.  Let $I$ and $J$ be subsets of
$\{1,\dots, N\}$ with $|I|= |J| = t$.  The subalgebra of $\omq$
generated by $x_{ij}$ with $i\in I$ and $j\in J$ can be regarded as an
algebra of $t\times t$ quantum matrices, and so we can calculate its
quantum determinant - this is a {\em $t\times t$ quantum minor} and we
denote it by $[I|J]$. 
The quantum minor $[I|I]$ is said to be a {\em principal quantum minor}.  
We denote the sum of all the principal
quantum minors of a given size $i$ by $\si$.  Note that $\sigma_1 =
x_{11} + \dots + x_{NN}$ and that $\sigma_N = \detq$.  
It is shown in \cite{dl} that the $\si$ are $\alpha$-coinvariants that
pairwise commute, and if $q$ is not a 
root of unity, then 
$\omq^{\coa}=\mc[\sigma_1,\ldots,\sigma_N]$, an $N$-variable commutative 
polynomial subalgebra of $\omq$. 
For the coaction $\beta$ we have that 
the {\em weighted sums of principal minors} 
$\ti:= \sum_{I}\, q^{-2w(I)}[I|I]$ (here 
$w(I)$ denotes the sum of the elements of $I$, 
and the summation ranges over all subsets $I$ of size $i$), 
$i=1,\ldots,N$, 
are pairwise commuting $\beta$-coinvariants. 
Moreover,  
assuming again that $q$ is not a root of unity,  
$\omq^{\cob}=\mc[\tau_1,\ldots,\tau_N]$, the $N$-variable commutative 
polynomial subalgebra generated by the $\ti$. 

We shall write $f(\xi)$ for $\evx(f)$, where $f\in\omq$ and 
$\xi$ is a ${\mathbb{C}}$-point of $M_q$. 

\begin{proposition}\label{ker-contains-coinvariants} 
The kernel of $\bx$ contains the right ideals  
\[\sum_{f\in\omq^{\cob}}(f-f(\xi))\omq\supseteq 
\sum_{i=1}^N(\ti-\ti(\xi))\omq.\] 
If $q$ is not a root of unity (or $q=1$), then the latter two right ideals 
are equal. 
The kernel of $\ax$ contains the left ideals  
\[\sum_{f\in\omq^{\coa}}\omq(f-f(\xi))\supseteq 
\sum_{i=1}^N\omq(\si-\si(\xi)).\] 
If $q$ is not a root of unity (or $q=1$), then the latter two 
left ideals are equal. 
\end{proposition} 

\begin{proof} If $f\in\omq^{\cob}$, then 
$\bx(f)=(\evx\otimes\id)(\beta(f))=(\evx\otimes\id)(f\otimes 1)
=f(\xi)\cdot 1\in\omq$, hence $f-f(\xi)\in\ker(\bx)$. By formula 
\eqref{beta-of-product} then $(f-f(\xi))\omq\subseteq \ker(\bx)$. 
Note that $\ti$ are $\beta$-coinvariants, and 
if $q$ is not a root of unity, they generate $\omq^{\cob}$ by \cite{dl}. 
Hence in the latter case the  
elements $\ti-\ti(\xi)$ generate the same right ideal within $\omq^{\cob}$ 
as all the $f-f(\xi)$ with $f\in\omq^{\cob}$. 

The same argument works for $\alpha$. 
\end{proof} 

\begin{remark}\label{kostant} 
{\rm More is known in the classical case $q=1$. If the adjoint 
orbit of $\xi\in\mnc$ 
is maximal, then $\ker(\bx)$, the vanishing ideal of the closure of the orbit 
of $\xi$, is generated by $\ti-\ti(\xi)$, $(i=1,\ldots,N)$; 
this is proved in \cite[Theorem 10]{k}, considering more generally 
the adjoint action on reductive Lie algebras. }
\end{remark} 

The coordinate Hopf algebra $\od$ of the diagonal subgroup $D$ of $\glnc$ 
is the commutative algebra  
$\mathbb{C}[t_1,t_1^{-1},\ldots,t_N,t_N^{-1}]$ of Laurent polynomials 
with comultiplication $\Delta(t_i)=t_i\otimes t_i$ and counit 
$\varepsilon(t_i)=1$. There exists a surjective Hopf algebra homomorphism 
$\pd:\ogq\to \mathcal{O}(D)$ determined by 
$\pd(x_{ij})=\delta_{ij}t_i$, $i,j=1,\ldots,N$. Therefore we say that 
$D$ is a quantum subgroup of $GL_q$, 
called the {\it diagonal subgroup}. 

There is a natural left coaction $\ld:=(\pd\otimes\id)\circ\Delta$ of 
$\od$ on $\ogq$.  The subset of $\ld$-coinvariants  is denoted by 
\[\odpergq=\{f\in\ogq\mid \ld(f)=1\otimes f\}.\] 
This is a subalgebra (since $\ogq$ is a comodule algebra with respect to 
$\ld$), as well as a right $\ogq$-subcomodule (since $\ld$ commutes 
with the right coaction of $\ogq$ on itself). So $\odpergq$ is  
{\it the coordinate algebra of a quantum homogeneous space} 
in the sense of \cite[11.6.1.]{ks}. 
This quantum homogeneous space is often called the
{\it quantum quotient space} $D\setminus G_q$, 
because it arises as the subalgebra of coinvariants with respect to the 
coaction of a quantum subgroup. 

Recall that in the classical case $q=1$, the subgroup $D$ is the stabilizer 
of $\xi$, provided that $\xi$ is a diagonal matrix with pairwise different 
eigenvalues, and the orbit of a diagonal $\xi$ is closed in 
$\omnc$. Hence for such $\xi$, the coordinate ring of the orbit of $\xi$ 
is $\mathcal{O}(D\setminus G_1)$. 
A partial analogue of this holds for any $q$. 
 
\begin{proposition}\label{image-in-diag-coinv} 
Let $\xi$ be a diagonal $N\times N$ matrix with complex entries. 
Then the image of $\bx$ is contained in $\odpergq$. 
Similarly we have $\im(\ax)\subseteq\odpergq$. 
\end{proposition} 

\begin{proof} 
The restriction of $\pd$ to $\omq$ is denoted by $\pc$, it maps 
$\omq$ onto $\oc=\mc[t_1,\ldots,t_N]$. By our assumption on $\xi$ the 
homomorphism $\evx:\omq\to\mc$ factors through $\pc$, we write 
$\evx$ also for the homomorphism $\oc\to\mc$, $t_i\mapsto \xi_{ii}$.  

By definition 
$\bx$ is the composition 
\begin{equation}\label{eq:beta-xi} 
\bx:\ \omq\stackrel{\beta}\longrightarrow\omq\otimes\ogq
\stackrel{\pc\otimes\id}\longrightarrow\oc\otimes\ogq
\stackrel{\evx\otimes\id}\longrightarrow\ogq,
\end{equation}  
and $\ld$ is the composition 
\begin{equation}\notag
\ld:\ \ogq\stackrel{\Delta}\longrightarrow\ogq\otimes\ogq
\stackrel{\pd\otimes\id}\longrightarrow\od\otimes\ogq.
\end{equation}
By linearity of $(\pd\otimes\id)\otimes\Delta$ we have 
\[(\pd\otimes\id)\circ\Delta\circ(\evx\otimes\id)
=(\evx\otimes\id\otimes\id)\circ(\id\otimes\pd\otimes\id)
\circ(\id\otimes\Delta)\] 
as maps 
$\oc\otimes\ogq\to\od\otimes\ogq$. 
Thus 
\[
(\pd\otimes\id)\circ\Delta\circ(\evx\otimes\id)\circ(\pc\otimes\id)=
(\evx\otimes\id\otimes\id)\circ(\pc\otimes\pd\otimes\id)\circ
(\id\otimes\Delta)\]  
as maps $\omq\otimes\ogq\to \od\otimes\ogq$.  
Therefore we have 
\begin{align*}
\ld\circ\bx 
=(\evx\otimes\id\otimes\id)\circ(\pc\otimes\pd\otimes\id)\circ
(\id\otimes\Delta)\circ\beta
\\=(\evx\otimes\id\otimes\id)\circ(\pc\otimes\pd\otimes\id)\circ
(\beta\otimes\id)\circ\beta  
\end{align*} 
(in the second equality we use the fact that $\omq$ is a right $\ogq$-comodule 
under $\beta$). 
Note that the formula \eqref{eq:beta} makes sense if we replace $\ogq$ 
by any Hopf algebra, in particular, it defines a right coaction 
$\beta^D:\od\to\od\otimes\od$. 
Since $\pd$ is a Hopf algebra homomorphism, we have 
$\beta^D\circ\pd=(\pd\otimes\pd)\circ\beta$. 
On the other hand, $\beta^D$ is equal to the trivial corepresentation 
$\id\otimes 1$, since it is 
the comorphism of the conjugation action of $D$ on 
itself, which is trivial. 
It follows that 
$(\pc\otimes\pd)\circ\beta=(\id\otimes 1)\circ\pc:\omq\to\oc\otimes\od$,  
and hence 
$\ld\circ\bx$ can be written as the composition 
\begin{eqnarray}\label{eq:lambda-circ-beta-xi}
\omq\stackrel{(\pc\otimes\id)\circ\beta}\longrightarrow\oc\otimes\ogq
\stackrel{(\id\otimes 1)\otimes\id}\longrightarrow
\\\notag   \oc\otimes\od\otimes\ogq
\stackrel{\evx\otimes\id\otimes\id}\longrightarrow \od\otimes\ogq. 
\end{eqnarray}

In order to prove our Proposition 
we need to show that for all $f\in\omq$ the equality 
$\ld(\bx(f))=1\otimes \bx(f)$ holds. 
Write $(\pd\otimes\id)(\beta(f))$ as $\sum_ia_i\otimes b_i$, 
where $a_i\in \oc$ and $b_i\in\ogq$. Then 
$\bx(f)=\sum_i\evx(a_i)b_i$ by \eqref{eq:beta-xi}. 
Therefore, using \eqref{eq:lambda-circ-beta-xi}, we have 
$\ld(\bx(f))=\sum_i\evx(a_i)1\otimes b_i=1\otimes\sum_i\evx(a_i)b_i
=1\otimes\bx(f)$. 
\end{proof} 


\section{Co-orbit maps of $2\times 2$ quantum matrices}\label{2x2} 

In the special case $N=2$ and $q$ not a root of unity (or $q=1$) 
we are able to refine 
Propositions~\ref{ker-contains-coinvariants} and ~\ref{image-in-diag-coinv}. 
Throughout this section we assume that $N=2$, 
$\omq$, $\ogq$ denote the coordinate rings of $2\times 2$ quantum matrices 
and quantum $GL(2)$, 
and $\od=\mc[t_1^{\pm 1},t_2^{\pm 1}]$. 
We will see in the proof of Lemma~\ref{D-coinv=K-coinv} that 
as an algebra $\odpergq$ is generated by $x_{11}x_{21}/\detq$, 
$x_{12}x_{21}/\detq$, $x_{12}x_{22}/\detq$. 
We write $\odpergq\lem$ for the subspace spanned by products of length at 
most $m$ in these generators. 

\begin{proposition} \label{2x2-im-ker}
Assume that $q$ is not a root of unity (or $q=1$), and take 
$\xi:=\left[\begin{matrix}\xi_1&0\\0&\xi_2\end{matrix}\right].$ 
\begin{itemize}
\item[(i)] 
If $\xi_1-q^{2k}\xi_2\neq 0$ for $k=0,1,2,\ldots$, then 
$\im(\bx)=\odpergq$. 
\item[(ii)] If $\xi_1-\xi_2,\xi_1-q^2\xi_2,\ldots,\xi_1-q^{2m-2}\xi_2$ 
are non-zero and $\xi_1-q^{2m}\xi_2=0$, then 
$\im(\bx)=\odpergq\lem$.  
\item[(iii)] If $\xi_1-q^{2k}\xi_2\neq 0$ for $k=0,1,2,\ldots$, 
then $\ker(\bx)=\sum_{i=1}^2(\ti-\ti(\xi))\omq$. 
\item[(iv)] If $\xi_1-\xi_2,\xi_1-q^2\xi_2,\ldots,\xi_1-q^{2m-2}\xi_2$ are 
non-zero and $\xi_1-q^{2m}\xi_2=0$, then $\ker(\bx)$ is the sum of 
$\sum_{i=1}^2(\ti-\ti(\xi))\omq$ and the $\ogq$-comodule generated by 
$\{x_{21}^k\mid k>m\}$. 
\item[(i')] 
If $\xi_1-q^{-2k}\xi_2\neq 0$ for $k=1,2,\ldots$, then $\im(\ax)=\odpergq$. 
\item[(ii')] If $\xi_1-q^{-2}\xi_2,\ldots,\xi_1-q^{-2m}\xi_2$ are nonzero 
and $\xi_1-q^{-2m-2}\xi_2=0$, then $\im(\ax)=\odpergq\lem$.  
\item[(iii')] If $\xi_1-q^{-2k}\xi_2\neq 0$ for $k=1,2,\ldots$, then 
$\ker(\ax)=\sum_{i=1}^2\omq(\si-\si(\xi))$. 
\item[(iv')] If $\xi_1-q^{-2}\xi_2,\ldots,\xi_1-q^{-2m}\xi_2$ are nonzero 
and $\xi_1-q^{-2m-2}\xi_2=0$, then $\ker(\ax)$ is the sum of 
$\sum_{i=1}^2\omq(\si-\si(\xi))$ and the $\ogq$-comodule generated by 
$\{x_{21}^k\mid k>m\}$.
\end{itemize} 
\end{proposition}

As an immediate corollary we obtain the following. 

\begin{theorem}\label{beta-main} 
Assume that $q$ is not a root of unity (or $q=1$) 
and take a diagonal complex matrix 
$\xi$ which is not a scalar multiple of any integral power of 
$\left[\begin{matrix}q^2&0\\0&1\end{matrix}\right]$. 
Then both $\bx$ and $\ax$ map $\omq$ surjectively onto 
$\odpergq$. 
In particular, 
$\im(\bx)=\im(\ax)$ is a quantum homogeneous space. 
Moreover, the kernel of $\bx$ is the right ideal of $\omq$ 
generated by the $\beta$-coinvariants 
$\tau_1-\tau_1(\xi)$, $\tau_2-\tau_2(\xi)$. 
The kernel of $\ax$ is the left ideal of $\omq$ generated by 
$\sigma_1-\sigma_1(\xi)$, $\sigma_2-\sigma_2(\xi)$. 
\end{theorem}

\begin{remark} 
{\rm In some sense the map $\bx$ (respectively $\ax$) 
is the best possible we can hope for. 
Suppose that 
$\varphi:\omq\to\odpergq$ 
is a map 
which is both a $\mc$-algebra homomorphism and intertwines between the 
$\ogq$-corepresentations $\beta$ and $\Delta$. 
Then $\ker\varphi$ is a completely prime two-sided ideal in $\omq$, since 
$\odpergq$ is a domain. 
We shall see below that the trivial $\ogq$-corepresentation appears 
with multiplicity $1$ in $\odpergq$. Therefore 
$\ker(\varphi)$ must contain $\tau_1-\lambda_1$ and $\tau_2-\lambda_2$ 
with some $\lambda_1,\lambda_2\in\mc$. 
Elementary commutator computations (exploiting the fact that $\tau_1$ 
is non-central) show that $\varphi(x_{12})=\varphi(x_{21})=0$, 
$\varphi(x_{11})\in\mc$, and $\varphi(x_{22})\in\mc$. In particular, 
$\im(\varphi)=\mc$ for any $\varphi$ compatible with both the algebra and 
comodule structures. }
\end{remark} 

It will be convenient for the proof of Proposition~\ref{2x2-im-ker} 
to pass from $\ogq$ to 
$\ohq$, the Hopf algebra quotient of $\ogq$ modulo the ideal 
generated by $\detq-1$. We keep the symbol $\Delta$ to denote 
the comultiplication in $\ohq$. Denote by $\pi$ the natural epimorphism 
$\pi:\ogq\to\ohq$, and 
$a:=\pi(x_{11})$, $b:=\pi(x_{12})$, $c:=\pi(x_{21})$, $d:=\pi(x_{22})$. 
The group $K$ of diagonal matrices in $SL(2,\mc)$ is a quantum subgroup 
of $\ohq$ with surjective Hopf algebra homomorphism 
$\pk:\ohq\to\ok=\mc[z^{\pm 1}]$ given by 
$\pk(a)=z$, $\pk(b)=0=\pk(c)$, $\pk(d)=z^{-1}$. 
There is a natural left coaction 
$\lk:=(\pk\otimes\id)\circ\Delta$ 
and a right coaction 
$\rk:=(\id\otimes\pk)\circ\Delta$ of $\ok$ on $\ohq$. 
The subalgebra of $\lk$-coinvariants is denoted by 
$\okperhq$. 

We shall use the corepresentation theory of $\ohq$, 
the material below can be found for example in \cite[4.2]{ks}. 
We assume that $q$ is not a root of unity (or $q=1$). 
Then $\ohq$ is cosemisimple, so any $\ohq$-comodule decomposes as the direct 
sum of irreducible subcomodules. 
The degree $n$ homogeneous component $P_n$ 
($n\in\mn:={\mathbb{N}}\cup\{0\}$) 
of the subalgebra 
$\mc\langle a,b\rangle$ of $\ohq$ is an $n+1$-dimensional right 
$\ohq$-subcomodule. Denote the corresponding corepresentation by 
$T_{n/2}$. 
Then $T_l:P_{2l}\to P_{2l}\otimes \ohq$, $l\in\frac 12\mn$, is a complete 
list of irreducible corepresentations of $\ohq$ up to isomorphism, see for 
example \cite[4.2.1]{ks}. 
We adopt the following notation: if 
$T:V\to V\otimes\ohq$ is a corepresentation of $\ohq$, then write 
$V[n]:=\{v\in V\mid (\id\otimes\pk)\circ T(v)=v\otimes z^n\}$, 
where $n\in\mz$. 
Since any $\ok$-comodule decomposes as a direct sum of one-dimensional 
$\ok$-subcomodules, 
the comodule $V$ decomposes as $\bigoplus_{n\in\mz}V[n]$. 
Assume that $V$ is finite dimensional. By the {\it character} 
$\chi(T)$ of $T$ 
we mean $\pk(f)$, where $f\in\ohq$ 
is the sum of the diagonal matrix elements of 
$T$ with repect to an arbitrarily chosen basis of $V$. 
In other words, 
$\chi(T)=\sum_{n\in\mz}\dim_{\mc}(V[n])z^n\in\ok$. 
(Note that our use of the word `character' slightly differs from that of 
\cite[11.2.2]{ks}, where $f$ is called the character of $T$, and $\chi(T)=\pk(f)$ the 
character of the $\ok$-corepresentation $(\id\otimes\pk)\circ T$.) 
For the irreducible $\ohq$-corepresentations $T_l$ we have that 
\[
\dim(P_{2l}[n])=\begin{cases} 
&1,\mbox{ if }n=2l,2l-2,2l-4,\ldots,-2l\\
&0,\mbox{ otherwise},\end{cases}
\]
so 
$\chi(T_l)=z^l+z^{l-2}+\cdots+z^{-l}$. 
In particular, the set $\{\chi(T_l)\mid l\in\frac{1}{2}\mn\}$ is linearly 
independent, and the multiplicity of $T_j$ as a direct summand of $T$ 
equals the coefficient of $\chi(T_j)$ in $\chi(T)$, expressed as a linear 
combination of the elements $\chi(T_l)$. 

The quantum quotient space $\os:=\okperhq$ is a well studied object, 
belonging to the one-parameter family of {\em quantum $2$-spheres} 
introduced in \cite{p}. The facts we need about $\os$ can be found in 
\cite[4.5]{ks}. As a subalgebra of $\ohq$, the quantum $2$-sphere 
$\os$ is generated by the elements 
\[\xmin:=(1+q^2)^{1/2}ac,\ \ 
\xnul:=1+(q+q^{-1})bc,\ \ 
\xplus:=(1+q^2)^{1/2}db.\] 
The right $\ohq$-corepresentation on $\os$ is isomorphic 
to $T_0\oplus T_1\oplus T_2\oplus\cdots$. 
The linear span of the products of length at most $n$ 
in the generators $\xmin,\xnul,\xplus$, 
denoted by $W^n$, 
is a right $\ohq$-subcomodule. The 
corepresentation on $W^n$ is isomorphic to 
$\bigoplus_{l=0}^n T_l$. Moreover, the simple subcomodule 
of $W^n$ corresponding to the summand $T_n$ is the subcomodule 
generated by $\xmin^n$, because 
$W^n[2n]=\mc\xmin^n$ and 
$W^n[d]=0$ if $d>2n$. 
In particular, as an 
$\ohq$-comodule $\os$ is generated by the powers of $\xmin$. 

The reason that we are able to switch from $\ogq$ to $\ohq$ is in the 
following lemma: 

\begin{lemma}\label{D-coinv=K-coinv} 
The restriction of $\pi$ maps $\odpergq$ isomorphically onto 
$\okperhq$. 
\end{lemma} 

\begin{proof} The coaction $\ld$ defines a $\mz\oplus\mz$-grading on 
$\ogq$: for $(m,n)\in \mz^2$  set 
\[\ogq_{(m,n)}:=\{f\in\ogq\mid \ld(f)=t_1^mt_2^n\otimes f\},\] 
then $\ogq=\bigoplus_{(m,n)\in\mz^2}\ogq_{(m,n)}$. 
By definition $\odpergq=\ogq_{(0,0)}$. 
For the algebra generators of $\ogq$ we have 
$x_{11},x_{12}\in\ogq_{(1,0)}$, $x_{21},x_{22}\in\ogq_{(0,1)}$, 
and $\detq^{-1}\in\ogq_{(-1,-1)}$. Therefore  
$\odpergq$ contains $x_{11}x_{21}/\detq$, 
$x_{12}x_{21}/\detq$, $x_{12}x_{22}/\detq$. 
These elements are mapped via $\pi$ onto a system of algebra generators 
of $\okperhq$. 

Since $\detq-1$ is not homogeneous with respect to this $\mz^2$-grading and 
$\ogq$ is a domain, we have that $\ker(\pi)=(\detq-1)\ogq$ does not 
contain any homogeneous elements. In particular, the kernel of 
$\pi$ is disjoint from $\odpergq$. 
\end{proof} 

The $\ogq$-coaction $\beta$ is replaced by the $\ohq$-coaction 
\[\psi:=(\id\otimes\pi)\circ\beta:\omq\to\omq\otimes\ohq,\]
and $\alpha$ is replaced by 
\[\varphi:=(\id\otimes\pi)\circ\alpha:\omq\to\omq\otimes\ohq.\] 
For a $\mc$-point $\xi$ in $M_q$ the maps $\psx,\phx:\omq\to\ohq$ are 
defined in an obvious way: 
$\psx=(\evx\otimes\id)\circ\psi=\pi\circ\bx$ and 
$\phx=(\evx\otimes\id)\circ\alpha=\pi\circ\ax$. 
Propositions~\ref{morphism-of-comodules} and \ref{image-in-diag-coinv} 
combined with Lemma~\ref{D-coinv=K-coinv} translate to the following. 
The map $\psx$ intertwines between the  
$\ohq$-corepresentations $\psi$ and $\Delta$, whereas 
$\phx$ intertwines between the $\ohq$-corepresentations 
$\varphi$ and $\Delta$. 
We have that $\ker(\psx)=\ker(\bx)$ and $\ker(\phx)=\ker(\ax)$.  
If $\xi$ is diagonal then $\im(\psx)\subseteq\okperhq$ 
with equality if and only if $\im(\bx)=\odpergq$, 
and $\im(\phx)\subseteq\okperhq$ with equality if and only if 
$\im(\ax)=\odpergq$. 
Moreover, 
$\im(\bx)=\odpergq^{\leq n}$ if and only if 
$\im(\psx)=W^n$, and 
$\im(\ax)=\odpergq^{\leq n}$ if and only if 
$\im(\phx)=W^n$. 

\bigskip
\begin{proofof}{Proposition~\ref{2x2-im-ker}} 

The kernel of 
$\psx$ contains $\sum_{i=1}^2(\ti-\ti(\xi))\omq$ 
by Proposition~\ref{ker-contains-coinvariants}, therefore  
$\psx$ factors through the composition of the natural homomorphisms 
\[\omq\stackrel{\mu}\longrightarrow C
\stackrel{\eta}\longrightarrow B,\] 
where 
$C:=\omq/(\tau_2-\tau_2(\xi))\omq$,  
$B:=C/uC$ 
with 
$u:=\mu(\tau_1-\tau_1(\xi))$. 
Set $y_{ij}:=\mu(x_{ij})$. 
Note that $\tau_2-\tau_2(\xi)=\detq-\xi_1\xi_2$ 
is a central element in $\omq$, so $C$ is a quotient algebra of $\omq$. 
It is a domain by \cite{j}. 
Since $\tau_2-\tau_2(\xi)$ is a $\psi$-coinvariant, 
$(\tau_2-\tau_2(\xi))\omq$ 
is a subcomodule of $\omq$ by 
\eqref{beta-of-product}, 
hence $C$ is a factorcomodule of $\omq$; 
the corresponding corepresentation is denoted by $\psi_C$. 
Similarly, $B$ is a factorcomodule of $C$, with corepresentation $\psi_B$. 
Our next aim is to show that 
\begin{equation}\label{eq:decomposition-of-B} 
\psi_B\cong T_0\oplus T_1\oplus T_2\oplus\cdots, 
\end{equation}
and moreover that, 
$\eta(y_{21})^r$ generates the simple subcomodule of $B$ corresponding to 
$T_r$ for $r\in\mn$. 

Denote by $\omq^{\leq r}$ 
the linear span of monomials of degree $\leq r$ in the 
generators $x_{ij}$, and 
$C^r:=\mu(\omq^{\leq r})$, 
$B^r:=\eta(C^r)$. 
Note that $\omq^{\leq r}$, $C^r$, $B^r$ are subcomodules. 
The symbols $\psi_{C^r}$, $\psi_{B^r}$ stand for the obvious 
subcorepresentations of $\psi_C$, $\psi_B$. 

For $r\in\mn$, a vector space basis of $C^r$ is 
\begin{equation}\label{eq:basis}
\{y_{11}^{i+1}y_{12}^jy_{21}^k,\ y_{12}^ly_{21}^my_{22}^n\mid 
i,j,k,l,m,n\in\mn,\ 
i+j+k\leq r-1,\ l+m+n\leq r\}  
\end{equation} 
(by the same argument as the corresponding result for $\ohq$ 
is proved in \cite[4.1.5]{ks}).  
From \eqref{eq:basis} we see that the character of the corepresentation 
$\psi_{C^r}$ is 
\begin{equation}\label{eq:character-of-C^r} 
\chi(\psi_{C^r})=
\sum_{i+j+k\leq r-1}z^{2(k-j)}+\sum_{l+m+n\leq r}z^{2(m-l)}.
\end{equation}  
(This character does not depend on $q$ or $\xi$.) 
In particular, 
$C^r[n]=0$ if $n>2r$, and 
$C^r[2r]=\mc y_{21}^r$, 
implying that the corepresentation 
$\psi_{C^r}$ has a unique irreducible summand isomorphic to $T_r$, 
and $y_{21}^r$ generates the corresponding simple subcomodule of $C^r$. 
It follows that $\eta(y_{21}^r)$ is either zero, or it generates a simple 
subcomodule of $B$, on which the given corepresentation is 
isomorphic to $T_r$. 

In order to show that $\eta(y_{21}^r)$ is non-zero, 
we introduce a $\mz$-grading $\deg$ on $\omq$ as follows: 
$\deg(x_{11}):=1$, $\deg(x_{22}):=-1$, 
$\deg(x_{12}):=0$, $\deg(x_{21}):=0$. 
It is easy to check that $\deg$ extends to an algebra-grading of $\omq$. 
Moreover, $\ker(\mu)$ is a homogeneous ideal, 
so $C$ inherits the grading. Note that $u$ is non-homogeneous, 
and $\deg(y_{21}^r)=0$. 
Since $C$ is a domain, any non-zero multiple of $u$ is non-homogeneous. 
Thus $y_{21}^r$ is not contained in $uC$. In other words, $\eta(y_{21}^r)$ 
is non-zero for any $r\in\mn$. 

Consequently, $\psi_{B^r}$ contains a subcorepresentation 
isomorphic to $T_0\oplus T_1\oplus\cdots T_r$, on the subcomodule generated by 
$1,\eta(y_{21}),\ldots,\eta(y_{21}^r)$. 
On the other hand, $\ker(\eta)\cap C^r=uC\cap C^r$ clearly contains 
$uC^{r-1}$, so $B^r$ is a factorcomodule of $C^r/uC^{r-1}$. 
Since $u$ is a coinvariant and $C$ is a domain, by 
\eqref{beta-of-product} we get that the corepresentation on $uC^{r-1}$ is 
isomorphic to $\psi_{C^{r-1}}$, hence the character of the corepresentation on 
$C^r/uC^{r-1}$ is $\chi(\psi_{C^r})-\chi(\psi_{C^{r-1}})$. 
Using \eqref{eq:character-of-C^r} we get that 
\begin{align}\label{eq:difference}
\chi(\psi_{C^r})-\chi(\psi_{C^{r-1}})
&=\sum_{i+j+k=r-1}z^{2(k-j)}+\sum_{l+m+n=r}z^{2(m-l)}
\notag\\&=\sum_{j+k\leq r-1}z^{2(k-j)}+\sum_{l+m\leq r}z^{2(m-l)}
\notag\\&=\sum_{s=0}^r\chi(T_s). 
\end{align}
So the corepresentation on $C^r/uC^{r-1}$ is also 
isomorphic to $\bigoplus_{s=0}^rT_r$. 
It follows that $C^r\cap uC=uC^{r-1}$, and 
$\psi_{B^r}\cong\bigoplus_{s=0}^r T_r$. 
Moreover, the subcomodule corresponding to $T_s$ is generated by 
$\eta(y_{21}^s)$. 
This shows \eqref{eq:decomposition-of-B}, since 
$B=\bigcup_{r=0}^{\infty}B^r$. 

As an immediate consequence of \eqref{eq:decomposition-of-B} 
and the fact that $\psx=\nu\circ\eta\circ\mu$ for some morphism 
$\nu:B\to \os$ of comodules,   
we obtain that the image of $\psx$ is the subcomodule of $\os$ generated 
by $\{\nu(\eta(y_{21}^n))\mid n\in\mn\}=\{\psx(x_{21}^n)\mid n\in\mn\}$.  
Moreover, the kernel of $\psx$ is $(\eta\circ\mu)^{-1}(\ker(\nu))$. 
By \eqref{eq:decomposition-of-B} 
$B$ is the direct sum of pairwise non-isomorphic simple subcomodules, 
and so any subcomodule is a direct sum of some of 
these simple subcomodules; in particular, 
$\ker(\nu)$ is the $\ohq$-subcomodule of $B$ generated by 
$\{\eta(y_{21}^n)\mid \nu(\eta(y_{21}^n))=0\}$. 
Therefore 
$(\eta\circ\mu)^{-1}(\ker(\nu))$ 
is the sum of 
$\ker(\eta\circ\mu)=\sum_{i=1,2}(\ti-\ti(\xi))\omq$ and 
the subcomodule of $\omq$ generated by 
$\{x_{21}^n\mid \psx(x_{21}^n)=0\}$. 
So (i), (ii), (iii), and (iv) follow from 
the fact that $W^n\leq\os$ 
is generated by $\{\xmin^r\mid r=0,1,\ldots,n\}$ as an $\ohq$-comodule, 
and from the evaluation of 
$\psx(x_{21}^n)$ given below. 

Since $\psx$ is not an algebra homomorphism, it is helpful to write it 
as the composition 
$\psx=p\circ\Psi$, 
where $\Psi$ is the map 
\[\Psi:=((S\circ\pi)\otimes\evx\otimes\pi)\circ\Delta^{(2)}: 
\omq\to\ohq\otimes\ohq,\]
and $p:\ohq\otimes\ohq\to\ohq$ is the multiplication map in $\ohq$. 
Introduce a new multiplication $*$ on $\ohq\otimes\ohq$ given by 
$(x_1\otimes y_1)*(x_2\otimes y_2)=x_2x_1\otimes y_1y_2$. 
Observe that $\Psi$ is an algebra homomorphism from $\omq$ to 
$(\ohq\otimes\ohq,*)$. 
We have that 
\begin{equation}\label{eq:psi}
\Psi\left[\begin{matrix}x_{11}&x_{12}\\x_{21}&x_{22}\end{matrix}\right]
=\left[\begin{matrix} 
\xi_1d\otimes a-q^{-1}\xi_2b\otimes c&\xi_1d\otimes b-q^{-1}\xi_2b\otimes d\\
-q\xi_1c\otimes a+\xi_2a\otimes c&-q\xi_1c\otimes b+\xi_2a\otimes d
\end{matrix}\right]. 
\end{equation}
By induction on $n$ we show that 
\begin{equation}\label{eq:c^n} 
\psx(x_{21}^n)=(-q)^n(\xi_1-\xi_2)(\xi_1-q^2\xi_2)\cdots 
(\xi_1-q^{2n-2}\xi_2)
c^na^n.
\end{equation} 
Indeed, \eqref{eq:c^n} holds for $n=1$ by \eqref{eq:psi}. 
Assume that \eqref{eq:c^n} holds for $n$. Then  
\begin{align*}\psx(x_{21}^{n+1})&=p(\Psi(x_{21}^{n+1}))
=p(\Psi(x_{21}^{n})*\Psi(x_{21}))
\\&=
p(\Psi(x_{21}^n)*(-q\xi_1c\otimes a+\xi_2a\otimes c))
\\&=
-q\xi_1cp(\Psi(x_{21}^{n}))a
+\xi_2ap(\Psi(x_{21}^n))c
\\&=(-q\xi_1+q^{2n+1}\xi_2)c\psx(x_{21}^n)a,
\end{align*} 
and formula \eqref{eq:c^n} easily follows for $n+1$. 
Hence if $\prod_{i=0}^{n-1}(\xi_1-q^{2i}\xi_2)\neq 0$, then 
$\psx(x_{21}^n)$ is a non-zero scalar multiple of $\xmin^n$, 
and if 
$\prod_{i=0}^{n-1}(\xi_1-q^{2i}\xi_2)=0$, then 
$\psx(x_{21}^n)=0$.

The same arguments show (i'), (ii'), (iii'), (iv'), by using the fact that 
\[\phx(x_{21}^n)=(-q)^n(\xi_1-q^{-2}\xi_2)(\xi_1-q^{-4}\xi_2)\cdots
(\xi_1-q^{-2n}\xi_2)a^nc^n.\] 

\end{proofof}


\section{Co-orbit maps for non-diagonal $2\times 2$ matrices} 

In this section we still restrict to the case $N=2$. 
The kernel and the image of the co-orbit map can be described also 
when $\xi$ is a non-diagonal $\mc$-point of $M_q$.  

\begin{proposition}\label{non-diagonal} 
Let $0\neq q\in\mc$, and take 
$\xi:=\left[\begin{matrix}0&\xi_1\\0&0\end{matrix}\right]$,  
where $\xi_1\neq 0$. 
\begin{itemize}
\item[(i)] The image of $\bx$ is the subalgebra of $\ogq$ generated by 
$x_{21}^2/\detq$, $x_{21}x_{22}/\detq$, $x_{22}^2/\detq$. 
In particular, $\im(\bx)$ is a quantum homogeneous space. 
\item[(ii)] The kernel of $\bx$ is the right ideal of 
$\omq$ generated by $\tau_1$ and $\tau_2$. 
\item[(i')] $\im(\ax)=\mc\langle x_{21}^2/\detq,x_{21}x_{22}/\detq, 
x_{22}^2/\detq\rangle$. 
\item[(ii')] $\ker(\ax)=\sum_{i=1}^2\omq\si$. 
\end{itemize} 
\end{proposition}

Again it is more convenient to switch from the $\ogq$-coactions to the 
corresponding $\ohq$-coactions. 

\begin{lemma}\label{restriction} 
The restriction of $\pi$ maps 
$\mc\langle x_{21}^2/\detq,x_{21}x_{22}/\detq, x_{22}^2/\detq\rangle$
isomorphically onto 
$\mc\langle c^2,cd,d^2\rangle$. 
\end{lemma} 

\begin{proof} The subalgebras 
$\mc\langle x_{21},x_{22}\rangle$ of $\ogq$ 
and $\mc\langle c,d\rangle$ of $\ohq$ 
are isomorphic to the coordinate ring of the quantum plane: 
the generators $q$-commute, and satisfy no more relations, as 
can be seen from standard basis arguments. 
The two algebras in our statement are isomorphic to 
the subalgebra of the coordinate 
ring of the quantum plane generated by its quadratic elements, 
or in other words, the sum of the even degree homogeneous components. 
\end{proof} 

\begin{proofof}{Proposition~\ref{non-diagonal}} 
As before, 
$\psi=(\id\otimes\pi)\circ\beta$ 
and 
$\psx=(\evx\otimes\id)\circ\psi$. The statements about $\bx$ 
translate by Lemma~\ref{restriction} to the assertions 
$\ker(\psx)=\sum_{i=1}^2\ti\omq$ and 
$\im(\psx)=\mc\langle c^2,cd,d^2\rangle$. 

As in the proof of Proposition~\ref{2x2-im-ker}, decompose $\psx$ 
as $p\circ\Psi$, where 
$\Psi=((S\circ\pi)\otimes\evx\otimes\pi)\circ\Delta^{(2)}$ 
and $p$ is the multiplication map in $\ohq$. By definition we have 
\begin{equation}\label{eq:psi2}
\Psi\left[\begin{matrix}x_{11}&x_{12}\\x_{21}&x_{22}\end{matrix}\right]
=\left[\begin{matrix} 
\xi_1d\otimes c&\xi_1d\otimes d\\
-q\xi_1c\otimes c&-q\xi_1c\otimes d
\end{matrix}\right]. 
\end{equation}
Since $\Psi$ is an algebra homomorphism from $\omq$ to 
$(\ohq\otimes\ohq,*)$, it is easy to see from the 
images of the generators that 
$\im(\Psi)$ is contained in 
$\sum_{k=0}^{\infty}
\mc\langle c,d\rangle_k\otimes
\mc\langle c,d\rangle_k$, 
where $\mc\langle c,d\rangle_k$ denotes the degree $k$ component of 
$\mc\langle c,d\rangle$. This space is mapped by $p$ into 
$W:=\mc\langle c^2,cd,d^2\rangle$, showing that 
$\im(\bx)\subseteq W$. 
We have 
\[\psx(x_{21}^n)
=p\circ\Psi(x_{21}^n)
=p(\Psi(x_{21})^n)
=(-1)^nq^n\xi^nc^{2n},\]
implying that $\im(\psx)$ contains $c^{2n}$ for all $n\in\mn$. 
Similarly, $\psx(x_{21}^ix_{12}^jx_{11}^k)$ 
is a non-zero scalar multiple of 
$(c^2)^i(d^2)^j(dc)^k$, hence 
we obtain (i).  

Set $C:=\omq/\tau_2\omq$ 
and $B:=C/\mu(\tau_1)C$, 
where $\mu:\omq\to C$ is the natural homomorphism. 
These 
are graded homomorphic images of $\omq$ (endowed with the usual 
grading). By Proposition~\ref{ker-contains-coinvariants} 
the linear map $\bx:\omq\to W$ factors through $B$, denote by 
$\eta:B\to W$ the induced surjective linear map. 
We shall show that $\eta$ is an isomorphism (this is clearly equivalent to 
the assertion $\ker(\psx)=\sum_{i=1}^2\tau_i\omq$). 
Define $\omq^{\leq r}$ as the degree $\leq r$ part of $\omq$, 
and write $C^r$,$B^r$, $W^r$ for the image of $\omq^{\leq r}$ in $C$, $B$, 
$W$. We know from \cite{j} that $C$ is a domain, and $C^r$ 
has the basis \eqref{eq:basis}. By homogeneity of $\mu(\tau_1)$ 
we have $C^r\cap \mu(\tau_1)C=\mu(\tau_1)C^{r-1}$. 
Hence 
$\dim(B^r)=\dim(C^r)-\dim(C^{r-1})=(r+1)^2$ 
(the latter equality follows for example from \eqref{eq:difference}). 
Since $\dim(W^r)=(r+1)^2$ as well, the map 
$\eta\big\vert_{B^r}:B^r\to W^r$ is an isomorphism. 
This holds for all $r\in\mn$, hence $\eta:B\to W$ is an isomorphism. 

The same arguments work for $\ax$. 
\end{proofof}

A similar result holds when 
$\xi=\left[\begin{matrix}0&0\\\xi_1&0\end{matrix}\right]$ 
is lower triangular. 

\begin{remark} 
{\rm $\mc\langle c^2,cd,d^2\rangle$ and $\os$ are isomorphic as 
right $\ohq$-comodules. However, they are not isomorphic 
quantum homogeneous spaces, since they are not isomorphic 
as algebras. Indeed, $\os$ is generated by $\xmin$ and 
$\xplus$ as a $\mc$-algebra with unity. 
Assume that $\mc\langle c^2,cd,d^2\rangle$ is also generated by two 
elements. We may suppose that the generators are contained in the maximal 
ideal $(c^2,cd,d^2)$. Then the images of the generators span the 
$3$-dimensional $\mc$-vector space $(c^2,cd,d^2)/(c^2,cd,d^2)^2$. 
This is a contradiction. 

It is not difficult to check that $\ohq$ is a free left 
$\mc\langle c^2,cd,d^2\rangle$-module. Therefore by 
\cite[Theorem 1]{t} the right coideal subalgebra 
$\mc\langle c^2,cd,d^2\rangle$ of $\ohq$ can be realized as 
the space of coinvariants with respect to the left coaction 
of an appropriate quotient left $\ohq$-module coalgebra of $\ohq$ 
(called a left coisotropic quantum subgroup of $\ohq$ in 
\cite{bcgst}). This quotient left $\ohq$-module coalgebra can be viewed 
as the ``stabilizer'' of $\xi$. 
}
\end{remark} 


\section{The generic case}\label{section-generic} 

We return to quantum matrices of arbitrary size $N$,  
but assume that $q$ is transcendental over the base field. 
It turns out that if the classical adjoint orbit of $\xi\in\mnc$ 
is of maximal dimension, then the kernel of the (quantized) co-orbit map 
coincides with the subset predicted by the classical theory, 
and is given in terms of coinvariants. 
For a diagonal $\xi$ with pairwise different eigenvalues, 
the image of the co-orbit map also reflects the classical picture. 

Let $A$ be an arbitrary commutative ring and $q$ a unit in $A$. 
We replace $\mc$ by $A$ in the constructions of Sections \ref{intro} 
and \ref{orbits}, and define the $A$-bialgebra $\oamq$,  
the $A$-Hopf algebras $\oagq$, $\oad$, 
the right coaction 
$\ba:\oamq\to\oamq\otimes\oagq$, 
the $\ba$-coinvariants $\tau_1,\ldots,\tau_N\in\oamq$. 
We say that an $N\times N$ matrix $\xi$ is an {\it $A$-point of} $M_q$, 
if the entries of $\xi$ are elements of $A$ satisfying the relations 
\eqref{eq:defining-relations}. 
Such a $\xi$ determines an $A$-algebra homomorphism 
$\evx:\oamq\to A$, $x_{ij}\mapsto \xi_{ij}$, and the co-orbit map 
$\bax:\oamq\to\oagq$ is defined as in Definition~\ref{co-orbit-map}. 
The $A$-algebra $\oamq$ is graded, the generators $x_{ij}$ having degree $1$. 
Write $\oamq^{\leq d}$ for the subspace spanned by the homogeneous components 
of degree $\leq d$. 

Let $K$ be a subfield of $\mc$, and let $q$ be an indeterminate over $K$. 
We shall apply the above constructions in the case when $A=\krq$, 
the field of rational functions in $q$, or when 
$A=\klq$, the ring of Laurent polynomials. 
When $(A,q)=(K,1)$, then $\oamq$ becomes $\okm$, the coordinate ring of 
the space of $N\times N$ matrices over $K$, and $\oagq$ becomes 
the coordinate ring $\okg$ of the general linear group $GL(N,K)$. 

Denote by $\eta:\klq\to K$ the 
$K$-algebra surjection mapping $q$ to $1$. 
The symbol $\eta$ will stand also for the induced surjections 
$M(N,\klq)\to M(N,K)$, 
$\oklqg\to\okg$. 

\begin{lemma}\label{lemma-specialization} 
The elements $f_1,\ldots,f_s\in\oklqg$ are $\krq$-linearly independent 
in $\okrqg$, if 
$\eta(f_1),\ldots,\eta(f_s)$ are $K$-linearly independent in 
$\okg$. 
\end{lemma} 

\begin{proof} 
Assume that $a_1f_1+\cdots+a_sf_s=0$ 
is a non-trivial relation with $a_i\in\krq$. 
Multiplying by an appropriate element of $K[q]$ we may achieve that 
$a_i\in \klq$ for all $i$. 
Since $\klq$ is a unique factorization domain, cancelling an appropriate 
power of $q-1$ we ensure that not all the $a_i$ are contained 
in the ideal $\langle q-1\rangle$. 
Then apply the map $\eta$ to get a non-trivial $K$-linear relation 
$\sum_{i=1}^s\eta(a_i)\eta(f_i)=0$. 
\end{proof} 

We fix a $\klq$-point $\xi$ of $M_q$ such that the centralizer of $\eta(\xi)$ 
in $M(N,K)$ has dimension $N$; that is, the adjoint orbit of $\eta(\xi)$ 
is of maximal dimension. 
For a fixed $d\in\mn$ set 
\[
X^d:=\okrqm^{\leq d}, 
\ \quad X_0^d:=\oklqm^{\leq d}, 
\ \quad X_1^d:=\okm^{\leq d},\] 
\[b^d:=\bx_{\krq}\big\vert_{X^d}, 
\ \quad b_0^d:=\bx_{\klq}\big\vert_{X_0^d}, 
\ \quad b_1^d:=\beta^{\eta(\xi)}_K\big\vert_{X_1^d}, \]
\[Y^d:=\im(b^d)\subset\okrqg, 
\ Y_0^d:=\im(b_0^d)\subset\oklqg, 
\ Y_1^d:=\im(b_1^d)\subset \okg. 
\]

Then we have a commutative diagram 
\begin{equation}\label{eq:diagram} \notag
\begin{array}{ccccc} 
X^d &\supset &X_0^d &\stackrel{\eta}\longrightarrow &X_1^d\\
&&&&\\
\downarrow{b^d} & &\downarrow{b_0^d} & &\downarrow{b_1^d}\\
&&&&\\
Y^d &\supset &Y_0^d &\stackrel{\eta}\longrightarrow &Y_1^d .
\end{array} 
\end{equation} 

Denote by $Z_1^d$ the kernel of $b_1^d$. 
It is the intersection of $X_1^d$ and the ideal of $\okm$ generated by 
$\eta(\ti-\ti(\xi))$, $i=1,\ldots,N$, by the classical result 
\cite[Theorem 10]{k} (the base field is $\mc$ in this paper; 
since $\beta^{\eta(\xi)}$ is obtained from 
$\beta^{\eta(\xi)}_K$ by extending scalars, and since 
the $\tau_i$ are defined over the field of rational numbers, 
the cited results of \cite{k} hold over any subfield $K$ of $\mc$). 
Moreover, by Lemma~\ref{lemma-comm-alg} below, 
$Z_1^d$ has a $K$-basis of the form 
\[\Lambda_1:=\bigcup_{i=1}^N\{\eta(\ti-\ti(\xi))w\mid w\in\Gamma^i\},\]
where $\Gamma^i$ is an appropriate set of monomials in the variables 
$x_{kl}$ of degree $\leq d-i$, for $i=1,\ldots,N$. 

\begin{lemma}\label{lemma-comm-alg} 
Let $R$ be a finitely generated commutative polynomial algebra over $K$, 
endowed with the usual grading. For $r\in R$ write 
$r=\bar r+\hat r$, where $\bar r$ is the highest degree homogeneous 
component of $r$. 
Let $u_1,\ldots,u_N$ be given elements of $R$, and assume that 
$\bar u_1,\ldots,\bar u_N$ is a regular sequence in $R$. 
Then for all $f$ in the ideal generated by 
$u_1,\ldots,u_N$, there exist elements $f_i\in R$, $i=1,\ldots,N$, 
such that 
$f=\sum_{i=1}^N u_if_i$, 
and $\deg(u_if_i)\leq\deg(f)$, $i=1,\ldots,N$. 
\end{lemma} 

\begin{proof} Apply induction on $N$. The case $N=1$ is trivial. 
Assume that the lemma is true for $N-1$. 
Write $f\in\sum_{i=1}^N u_iR$ as 
$f=\sum_{i=1}^Nu_ih_i$, 
where $d:=\mathrm{max}\{\deg(u_ih_i)\mid i=1,\ldots,N\}$ 
is minimal. 
If $d\leq \deg(f)$, then we are done. 
Suppose that $d>\deg(f)$. 
Without loss of generality we may assume that $\deg(u_Nh_N)=d$ 
because any permutation of $\bar u_1,\ldots,\bar u_N$ is a regular sequence, 
see the corollary to \cite[Theorem 16.3]{ma}. 
(We note that it is not essential for the proof 
to make this assumption, however, it simplifies the notation below.) 
Then 
\[\sum_{i:\deg(u_ih_i)=\deg(u_Nh_N)}\bar u_i\bar h_i=0.\]  
By assumption $\bar u_N$ is not a zero-divisor modulo 
$\sum_{i=1}^{N-1}\bar u_iR$, hence 
$\bar h_N=\sum_{i=1}^{N-1}\bar u_ig_i$ with appropriate 
homogeneous elements $g_i$, where 
$\deg(\bar u_ig_i)=\deg(\bar h_N)$, $i=1,\ldots,N-1$. 
Thus 
\begin{align*}
u_Nh_N &=u_N(\bar h_N+\hat h_N)
\\&=u_N\hat h_N+\sum_{i=1}^{N-1}u_N\bar u_ig_i
\\&=u_N\hat h_N+\sum_{i=1}^{N-1}u_N(u_i-\hat u_i)g_i
\\&=\sum_{i=1}^{N-1}u_i(u_Ng_i)+u_N(\hat h_N-\sum_{i=1}^{N-1}\hat h_ig_i). 
\end{align*} 
Obviously $\deg(u_iu_Ng_i)\leq d$, 
$\deg(u_N\hat u_ig_i)<d$ 
for $i=1,\ldots,N-1$. 
Set 
$b_i:=h_i+u_Ng_i$ for $i=1,\ldots,N-1$, and 
$b_N:=\hat h_N-\sum_{i=1}^{N-1}\hat u_ig_i$. 
We have 
$f=\sum_{i=1}^Nu_ib_i$, 
with $\deg(u_ib_i)\leq d$ for $i=1,\ldots,N-1$, 
and 
$\deg(u_Nb_N)<\deg(u_Nh_N)=d$. 
We claim that 
$\deg(\sum_{i=1}^{N-1}u_ib_i)<d$. 
Indeed, either 
$\deg(\sum_{i=1}^{N-1}u_ib_i)\leq\deg(u_Nb_N)<d$, 
or 
$\deg(\sum_{i=1}^{N-1}u_ib_i)>\deg(u_Nb_N)$ 
implying 
$\deg(\sum_{i=1}^{N-1}u_ib_i)=\deg(f)<d$. 
By the induction hypothesis on $N$, there exist 
$c_1,\ldots,c_{N-1}\in R$ with 
$\sum_{i=1}^{N-1}u_ib_i=\sum_{i=1}^{N-1}u_ic_i$, 
and 
$\deg(u_ic_i)\leq\deg(\sum_{i=1}^{N-1}u_ib_i)<d$ 
for $i=1,\ldots,N-1$. 
Therefore 
\[f=u_1c_1+\cdots+u_{N-1}c_{N-1}+u_Nb_N,\] 
where each summand on the right hand side has degree $<d$. 
This contradicts to our assumption on the minimality of $d$. 
\end{proof} 

Since $\eta(\tau_1),\ldots,\eta(\tau_N)$ is a regular sequence in 
$\okm$ by \cite[Theorem 10]{k}, 
we may apply Lemma~\ref{lemma-comm-alg} to conclude the existence 
of the basis $\Lambda_1$ in $Z_1^d$. 
Our next aim is to show that $\Lambda_1$ can be lifted to a basis of 
the kernel of $b^d$. 
There is a $\mz^N$-grading on the spaces 
$X^d$, $X_1^d$, $Y^d$, $Y_1^d$, 
determined by a right coaction of $\oad$. 
This coaction is 
$(\id\otimes\pd)\circ\ba$ for $X^d$, $X_1^d$, 
and 
$(\id\otimes\pd)\circ\Delta$ for $Y^d$, $Y_1^d$. 
(Recall that a coaction $\varphi:V\to V\otimes\oad$ 
yields the direct sum decomposition 
$V=\bigoplus_{a\in\mz^N}\{v\in V\mid 
\varphi(v)=v\otimes t_1^{a_1}\cdots t_N^{a_N}\}$.)  
The maps $b^d$, $b_1^d$, $\eta$ are compatible with the $\mz^N$-grading, 
because $b^d$, $b_1^d$ are homomorphisms of $\oagq$-comodules, hence they are 
homomorphisms of $\oad$-comodules, whereas $\eta$ maps a monomial in 
$\oklqg$ to formally the same monomial in $\okg$. 
Write $\chi(V)$ for the Hilbert series of a finite dimensional 
$\mz^N$-graded vector space $V$, so 
$\chi(V)$ is an element of $\mz[t_1^{\pm 1},\ldots,t_N^{\pm 1}]$. 
In other words, $\chi(V)$ is the character of the 
corresponding corepresentation of $\oad$ on $V$. 
Introduce a partial order $\geq$ on 
$\mz[t_1^{\pm 1},\ldots,t_n^{\pm 1}]$ 
by defining $f\geq g$ if all coefficients of $f-g$ are non-negative. 
Obviously $W\leq V$ and $\chi(W)=\chi(V)$ imply $W=V$. 

Denote by $Z^d$ the intersection of $X^d$ and the right ideal of 
$\okrqg$ generated by $\ti-\ti(\xi)$, $i=1,\ldots,N$. 
The same argument as in the proof of 
Proposition~\ref{ker-contains-coinvariants} 
shows that $Z^d\subseteq \ker(b^d)$. 

\begin{proposition}\label{prop-hilbert-series} 
For all $d\in\mn$ we have the following:
\begin{itemize}
\item[(i)] the kernel of $b^d$ is $Z^d$; 
\item[(ii)] the Hilbert series $\chi(Y^d)$ equals $\chi(Y_1^d)$. 
\end{itemize} 
\end{proposition} 

\begin{proof} 
Choose a set $\Omega_1$ of monomials of degree $\leq d$ in $\okm$ 
such that $\Lambda_1\cup\Omega_1$ is a $K$-basis of $X_1^d$. 
Set 
\[\Lambda:=\{(\ti-\ti(\xi))\hat w\mid i=1,\ldots,N;\ w\in\Gamma^i\},\] 
where $\hat w$ is a chosen monomial in $X_0^d$ with $\eta(\hat w)=w$. 
The elements of $\Lambda$ are multihomogeneous with respect to the 
$\mz^N$-grading introduced on $X^d$. 
Similarly, lift each $v\in\Omega_1$ to a monomial $\hat v\in X_0^d$ with 
$\eta(\hat v)=v$. We obtain the set 
$\Omega:=\{\hat v\mid v\in\Omega_1\}$ 
of multihomogeneous elements in $X^d$. 

As we noted above, $\Lambda$ is a subset of $\ker(b^d)$, and it is 
$\krq$-linearly independent by Lemma~\ref{lemma-specialization}. 
Therefore 
\begin{equation}\label{eq:*} 
\chi(\ker(b^d))\geq \chi(Z^d)\geq
\chi(\mathrm{Span}_{\krq}\{\Lambda\})=
\chi(\mathrm{Span}_K\{\Lambda_1\})
=\chi(Z_1^d).
\end{equation} 
Since $Z_1^d=\ker(b_1^d)$, we have that $\Omega_1$ is mapped 
under $b_1^d$ to a basis of $Y_1^d$. 
Again by Lemma~\ref{lemma-specialization} 
we have that $b^d(\Omega)$ is a $\krq$-linearly independent subset of $Y^d$. 
Thus 
\begin{equation}\label{eq:**} 
\chi(Y^d)\geq \chi(\mathrm{Span}_{\krq}\{b^d(\Omega)\})
=\chi(\mbox{Span}_K\{b_1^d(\Omega_1)\})=\chi(Y_1^d). 
\end{equation} 
It is well known (see \cite[9.2.1 Proposition 6]{ks}) 
that there is a set of monomials in the variables $x_{ij}$  
which is a basis both for $X^d$ and $X_1^d$. 
It follows that 
\begin{equation}\label{eq:***}
\chi(Y^d)+\chi(\ker(b^d))=\chi(X^d)
=\chi(X_1^d)=\chi(Y_1^d)+\chi(Z_1^d).
\end{equation} 
Comparing \eqref{eq:*}, \eqref{eq:**}, \eqref{eq:***} we obtain that 
all inequalities in \eqref{eq:*} and \eqref{eq:**} 
must be equalities. 
In particular, 
$\chi(Y^d)=\chi(Y_1^d)$ in \eqref{eq:**} 
and 
$\chi(\ker(b^d))=\chi(\mathrm{Span}_{\krq}\{\Lambda\})$ in 
\eqref{eq:*}, 
implying 
$\ker(b^d)=Z^d=\mathrm{Span}_{\krq}\{\Lambda\}$. 
\end{proof} 

As an immediate corollary we obtain the following. 

\begin{theorem}\label{thm-kernel} 
Let $\xi$ be a $\krq$-point of $M_q$ such that the centralizer of $\eta(\xi)$ 
in $M(N,K)$ has dimension $N$. 
Then the kernel of $\bx_{\krq}$ coincides with the right ideal 
$\sum_{i=1}^N(\ti-\ti(\xi))\okrqm$. 
\end{theorem} 

We determine the comodule structure of the image of $\bx_{\krq}$ 
with the aid of corepresentation theory. 
The material we shall summarize below can be found in \cite[11.5]{ks}
(the base field is $\mc$ in that book; however, the irreducible 
corepresentations of $\ogq$ are defined over the field of rational numbers, 
so the results mentioned below obviously hold over $K$ and $\krq$). 
In the sequel $(A,q)$ will stand for either of 
$(\krq,q)$ or $(K,1)$. 
Recall that $\oagq$ is cosemisimple. 
The irreducible corepresentations are indexed by the set $P$ 
of dominant integral weights for $GL(N)$. 
Denote by $T_q(\lambda)$ 
the irreducible $\oagq$-corepresentation 
associated with $\lambda\in P$, write $\chi(T_q(\lambda))$ for its character. 
Recall that by the {\it character} $\chi(T)$ of a finite dimensional 
$\oagq$-corepresentation $T:V\to V\otimes\oagq$ 
we mean the image under $\pd$ of the sum 
of the diagonal matrix coefficients. 
In other words, the character $\chi(T)$ 
is the same as the Hilbert series $\chi(V)$ 
with respect to the $\mz^N$-grading of $V$ 
determined by the $\oad$-corepresentation $(\id\otimes\pd)\circ T$. 
Any $\oagq$-corepresentation decomposes as a direct sum of copies of 
the irreducible corepresentations $T_q(\lambda)$. 
Given a finite dimensional corepresentation $T$, the multiplicity of 
$T_q(\lambda)$ as a summand of $T$ is the same as the 
(uniquely determined) coefficient of 
$\chi(T_q(\lambda))$ in $\chi(T)$, 
expressed as a linear combination of the characters 
$\chi(T_q(\mu))$, $\mu\in P$. 
We shall need also the fact that the character of 
the irreducible $\okrqg$-corepresentation $T_q(\lambda)$ coincides with 
the character of the $\okg$-corepresentation $T_1(\lambda)$ 
(which coincides with the character of the irreducible rational 
representation of $GL(N,\mc)$ associated with $\lambda$). 
Consequently, if we are given an $\okrqg$-corepresentation 
$T:V\to V\otimes \okrqg$ and an $\okg$-corepresentation 
$U:W\to W\otimes\okg$ such that $\chi(V)=\chi(W)$, then 
the multiplicity of $T_q(\lambda)$ as a summand of $T$ 
is the same as the multiplicity of $T_1(\lambda)$ as a summand of 
$U$, for all $\lambda\in P$. 

To simplify the notation we write $T(\lambda):=T_1(\lambda)$.  
The $\okg$-corepresentation on 
$\im(\beta^{\eta(\xi)}_K)$ 
decomposes as 
\[\im(\beta^{\eta(\xi)}_K)
\cong \bigoplus_{\lambda\in P}m(\lambda)T(\lambda),\] 
where $m(\lambda)$ denotes the dimension of the zero weight space in the dual 
corepresentation $T^*(\lambda)$,  
by \cite[Theorem 0.4 and formula (0.1.6)]{k}.  
Although the result is stated in \cite{k} in the 
language of representations, it can clearly be converted to the 
language of corepresentations. Indeed, observe that 
the character of the $\okg$-corepresentation on 
$\im(\beta^{\eta(\xi)}_K)$ coincides with the usual
formal character of the natural right action of $GL(N,K)$ on this space, 
and $\im(\beta^{\eta(\xi)}_K)$ is the coordinate ring of the 
closure of the adjoint orbit of $\eta(\xi)$. 
So $m(\lambda)$ is the dimension of the subspace of the underlying 
vector space of $T^*(\lambda)$ consisting of the vectors $v$ with 
$(\id\otimes\pd)\circ T^*(\lambda)(v)=v\otimes 1$. 
In particular, $m(\lambda)$ is finite for all $\lambda\in P$, and does not 
depend on $\eta(\xi)$ (still assuming that the centralizer of $\eta(\xi)$ 
has dimension $N$). 
The multiplicity of $T_q(\lambda)$ in the $\okrqg$-corepresentation on 
$\okrqg$ is the dimension of $T_q(\lambda)$ 
(see \cite[11.5.4 Theorem 51 and 11.1.4 Proposition 8 (ii)]{ks}), 
hence the subcorepresentation on 
$\im(\bx_{\krq})$ decomposes as 
\[\im(\bx_{\krq})\cong 
\bigoplus_{\lambda\in P} m^{\xi}_q(\lambda) T_q(\lambda)\]  
with $m^{\xi}_q(\lambda)$ finite for all $\lambda\in P$. 

\begin{proposition}\label{prop-multiplicities} 
Let $\xi$ be a $\krq$-point of $M_q$ such that the centralizer of 
$\eta(\xi)$ in $M(N,K)$ has dimension $N$. 
Then for all $\lambda\in P$ we have 
$m^{\xi}_q(\lambda)=m(\lambda)$. 
\end{proposition} 

\begin{proof} 
Denote by 
$m^d(\lambda)$ the multiplicity of $T(\lambda)$ as a summand in $Y_1^d$, 
and denote by 
$m^{\xi,d}_q(\lambda)$ the multiplicity of $T_q(\lambda)$ 
as a summand in $Y^d$. 
Obviously 
\begin{equation}\label{eq:sup1}
m(\lambda)=\sup\{m^d(\lambda)\mid d\in\mn\}
\end{equation} 
and 
\begin{equation} \label{eq:sup2} 
m^{\xi}_q(\lambda)=\sup\{m^{\xi,d}_q(\lambda)\mid d\in\mn\}, 
\end{equation} 
because 
$\im(\bx_K)=\bigcup_{d\in\mn} Y_1^d$ and 
$\im(\bx_{\krq})=\bigcup_{d\in\mn} Y^d$. 
On the other hand, 
$\chi(Y^d)=\chi(Y_1^d)$ by Proposition~\ref{prop-hilbert-series}, 
implying 
\begin{equation}\label{eq:mult=}
m^{\xi,d}_q(\lambda)=m^d(\lambda)
\mbox{ for all }d\in\mn,\  \lambda\in P.
\end{equation} 
So our statement follows from equations \eqref{eq:sup1}, 
\eqref{eq:sup2}, \eqref{eq:mult=}. 
\end{proof}

If $\eta(\xi)$ is diagonal with pairwise different eigenvalues, 
then the stabilizer of $\eta(\xi)$ 
with respect to the classical adjoint action is the diagonal 
subgroup of $GL(N,K)$, and the orbit of $\eta(\xi)$ is closed in $M(N,K)$. 
It follows then that the image of $\beta^{\eta(\xi)}_K$ is $\okdperg$. 
In particular, the multiplicity of $T(\lambda)$ as a direct summand 
of the $\okg$-corepresentation on $\okdperg$ is $m(\lambda)$. 
We show that the same holds in the $q$-deformed case. 

\begin{lemma}\label{lemma-mult-in-odpergq}
The multiplicity of $T_q(\lambda)$ as a direct summand of 
the $\okrqg$-corepresent\-ation on $\okrqdperg$ is $m(\lambda)$ for all 
$\lambda\in P$. 
\end{lemma} 

\begin{proof} 
In the sequel $(A,q)$ will stand for either of 
$(\krq,q)$ or $(K,1)$. 
First we describe a spanning set of $\oadperg$. 
For an $N\times N$ matrix $a=(a_{ij})$ with non-negative integer 
entries and a non-negative integer $d$ set 
\[x^a_d:=\detq^{-d}x_{11}^{a_{11}}x_{12}^{a_{12}}\cdots x_{NN}^{a_{NN}}\] 
(the variables $x_{ij}$ are ordered lexicographically). 
The elements of this form span $\oagq$, and 
\[(\pd\otimes\id)\circ\Delta(x^a_d)=t^a_d\otimes x^a_d,\] 
where 
\[t^a_d:=\prod_{i=1}^Nt_i^{-d+\sum_{j=1}^Na_{ij}}.\] 
It follows that $\oadperg$ is spanned by the $x^a_d$ with 
$\sum_{j=1}^Na_{ij}=d$ for $i=1,\ldots,N$. 
Set 
\[L^d_A:=\mathrm{Span}_A\{x^a_d\mid \sum_{j=1}^Na_{ij}=d,\ i=1,\ldots,N\}.\]
Observe that the given spanning set of $L^d_A$ is actually a basis, 
because after multiplication by $\detq^d$, 
it becomes a set of linearly independent 
monomials of degree $dN$ from the standard basis of the domain $\oamq$. 
Clearly 
\[\oadperg=\bigcup_{d\in\mn}L^d_A,\] 
and $L^d_A$ is a right $\oagq$-subcomodule of $\oadperg$ for all $d$. 

Note that 
\[(\id\otimes\pd)\circ\Delta(x^a_d)
=x^a_d\otimes t_{a,d}
\mbox{ \ with \ } 
t_{a,d}:=\prod_{j=1}^Nt_j^{-d+\sum_{i=1}^Na_{ij}}.\] 
It follows that 
\[\chi(L^d_A)=\sum_{a\in\omega_d} t_{a,d},\] 
where 
\[\omega_d:=\{a=(a_{ij})\in M(N,\mn)\mid 
\sum_{j=1}^Na_{ij}=d\mbox{ for }i=1,\ldots,N\}.\] 
Observe that the above computation of 
$\chi(L^d_A)$ is the same when $A=\krq$ and when $A=K$. 
It follows that for all $d\in\mn$ the multiplicity of 
$T_q(\lambda)$ in $L_{\krq}^d$ is the same as the 
multiplicity of $T(\lambda)$ in $L_K^d$. 
Consequently,  
the multiplicity of $T_q(\lambda)$ in 
$\okrqdperg$ is the same as the multiplicity of 
$T(\lambda)$ in $\okdperg=\im(\beta^{\eta(\xi)}_K)$. 
\end{proof} 

\begin{theorem}\label{thm-image} 
Let $\xi$ be a diagonal $\klq$-point of $M_q$ such that 
$\eta(\xi)$ has pairwise different eigenvalues. 
Then the image of $\bx_{\krq}$ coincides with 
$\okrqdperg$. 
\end{theorem} 

\begin{proof} 
We have 
$\im(\bx_{\krq})\subseteq \okrqdperg$ by the same argument as in the proof of 
Proposition~\ref{image-in-diag-coinv}. 
The multiplicity of $T_q(\lambda)$ is $m(\lambda)$ 
both in $\okrqdperg$ and in $\im(\bx_{\krq})$ by 
Proposition~\ref{prop-multiplicities} and Lemma~\ref{lemma-mult-in-odpergq}. 
This implies that $\im(\bx_{\krq})$ is the whole of 
$\okrqdperg$. 
\end{proof} 

Finally  we return to the coordinate ring of quantum matrices over $\mc$. 
The results of this section have the following corollary 
in the framework of Section~\ref{orbits}. 

\begin{corollary} 
Let $\xi$ be a $\mc$-point of $M_q$ such that the centralizer of 
$\xi$ in $M(N,K)$ has dimension $N$. 
Assume that $q\in\mc^*$ is transcendental over the subfield $K$ of $\mc$ 
generated by the entries of $\xi$. 
Then the kernel of 
$\bx$ coincides with the right ideal of $\omq$ 
generated by $\ti-\ti(\xi)$, $i=1,\ldots,N$. 
If in addition $\xi$ is diagonal, 
then the image of $\bx$ is $\odpergq$. 
\end{corollary}  

\begin{proof} 
Since $q$ is transcendental over $K$, we 
may apply Theorems~\ref{thm-kernel} and ~\ref{thm-image} 
for the kernel and the image of $\bx_{\krq}$. 
The $\mc$-linear map $\beta^{\xi}$ is obtained 
from the $\krq$-linear map $\bx_{\krq}$ by extending scalars to $\mc$. 
Therefore a generating set of $\ker(\bx_{\krq})$ 
as a $\krq$-vector space 
spans $\ker(\bx)$ over $\mc$, 
and a generating set of $\im(\bx)$ as a $\krq$-vector space 
spans $\im(\bx)$ over $\mc$. 
\end{proof} 

It is straightforward to modify the statements of this section 
to obtain similar results for the coaction $\alpha$. 


\noindent M. Domokos: 
R\'enyi Institute of Mathematics, Hungarian Academy of Sciences,\\ 
P.O. Box 127, 1364 Budapest, Hungary\\
E-mail: domokos@renyi.hu, domokos@maths.ed.ac.uk\\
\noindent(Domokos is in Edinburgh until February 2004.)
\\
\\
\noindent R. Fioresi: 
Dipartimento di Matematica, Universita' di Bologna\\
Piazza Porta San Donato 5, 40126 Bologna, Italy 
\\E-mail: fioresi@dm.unibo.it
\\
\\
\noindent T. H. Lenagan: 
School of Mathematics, University of Edinburgh,
\\ James Clerk Maxwell Building, King's Buildings, Mayfield Road, 
\\Edinburgh EH9 3JZ, Scotland
\\E-mail: tom@maths.ed.ac.uk

\end{document}